\documentclass{amsart}

\usepackage[latin1]{inputenc}
\usepackage[T1]{fontenc}
\usepackage{amsmath,amssymb,amsthm}
\usepackage{graphicx}
\usepackage{epsfig}
\usepackage{hyperref}
\usepackage[abs]{overpic}

\begin{document}

\newtheorem{theorem}{Theorem}[section]
\newtheorem{lemma}[theorem]{Lemma}
\newtheorem{corollary}[theorem]{Corollary}
\newtheorem{cor}[theorem]{Corollary}
\newtheorem{proposition}[theorem]{Proposition}
\theoremstyle{definition}
\newtheorem{definition}[theorem]{Definition}
\newtheorem{example}[theorem]{Example}
\newtheorem{claim}[theorem]{Claim}
\newtheorem{xca}[theorem]{Exercise}

\theoremstyle{remark}
\newtheorem{remark}[theorem]{Remark}

\newcommand{\pint}{P_{\text{int}}}
\newcommand{\B}{\partial(D)}
\newcommand{\Z}{\mathbb Z}
\newcommand{\N}{\mathbb N}
\newcommand{\I}{I_{n/2}}
\newcommand{\Q}{\mathbb Q}
\newcommand{\Ztwo}{\mathbb Z_2}
\newcommand{\T}{T}
\newcommand{\BB}{\partial(\T)}
\newcommand{\Btilde}{\partial(\tilde D)}
\newcommand{\BBtilde}{\partial(\tilde \BB)}
\newcommand{\Bprime}{\partial(D')}
\newcommand{\len}{L}
\newcommand{\diam}{\text{diam}}
\newcommand{\be}{\begin{equation}}
\newcommand{\ee}{\end{equation}}
\newcommand{\prob}{\mbox{\bf P}}
\newcommand{\vint}{V_{\text{int}}}
\newcommand{\dint}{D_{\text{int}}}
\newcommand{\tint}{T^{\text{int}}}
\newcommand{\fdis}{F_{\text{dis}}}
\newcommand{\vdis}{V_{\text{dis}}}
\newcommand{\fsid}{F_{\text{dis}}}
\newcommand{\vsid}{V_{\text{dis}}}
\newcommand{\basin}{{\text{Basin}}}
\newcommand{\Star}{{\text{Star}}}
\newcommand{\atc}{associated triangulated disk}
\newcommand{\ptc}{triangulated disk}

\newcommand\HFILL{\hspace*{\fill}}
\newcommand\VFILL{\vspace*{\fill}}

\newenvironment{pf} {\noindent{\sc Proof. }}{{\hfill
$\Box$}\par\vskip2\parsep}

\newenvironment{pfofthm}[1]
{\par\vskip2\parsep\noindent{\sc Proof of Theorem\ #1. }}{{\hfill
$\Box$}
\par\vskip2\parsep}

\newenvironment{pfoflem}[1]
{\par\vskip2\parsep\noindent{\sc Proof of Lemma\ #1. }}{{\hfill
$\Box$}
\par\vskip2\parsep}

\newenvironment{skpf} {\noindent{\sc Sketch
Proof. }}{{\hfill $\Box$}\par\vskip2\parsep}
\newenvironment{skpfof}[1] {\par\vskip2\parsep\noindent{\sc Sketch
Proof of\ #1. }}{{\hfill $\Box$}\par\vskip2\parsep}

\newcommand{\blankbox}[2]{%
  \parbox{\columnwidth}{\centering
    \setlength{\fboxsep}{0pt}%
    \fbox{\raisebox{0pt}[#2]{\hspace{#1}}}%
  }%
}

\bibliographystyle{plain}

\title{Exponential clogging time for a one dimensional DLA}

\author{Itai Benjamini and Christopher Hoffman}
\address{Department of Mathematics, The Weizmann Institute, Rehovot, Israel 76100 and Microsoft Research}
\email{itai.benjamini@weizmann.ac.il}
\address{Department of Mathematics, University of Washington, Seattle, WA 98195}
\email{hoffman@math.washington.edu}

\date{\today.}

\maketitle

\vspace*{-0.5in}

\begin{abstract}
In this paper a  simple DLA type model is analyzed. In \cite{BY}
the standard DLA model from \cite{WS} was considered on a cylinder
and the arm growing phenomena was established, provided that the
section of the cylinder has sufficiently fast mixing rate. When
considering DLA on a cylinder it is natural to ask how many
particles it takes to clog the cylinder, e.g. modeling clogging of
arteries. In this note we formulate a very simple DLA clogging
model and establish an exponential lower bound on the number of
particles arriving before clogging appears. In particular we
possibly shed some light on why it takes so long to reach the
bypass operation.\end{abstract}

\section{Introduction}

We start with an informal description of our model. Fix some $N
\in \N$. Initially there is a particle only at the vertices  $0
\mbox{ and } 1 \in \Z$. A third particle performs a simple random
walk started at positive infinity until at some random time when
it stops and never moves again. If the particle is at a vertex
$i$, it "freezes" there (stays there and remains there for all
time) with probability equal to the number of particles at vertex
$i-1$ divided by $N$. If the particle does not freeze at $i$ then
it takes a step in simple random walk.  The particle repeats this
procedure of either freezing or taking one step of simple random
walk until it freezes. When the third particle freezes, the forth
particle starts and so on. By the analysis in \cite{BY} the
cluster of the particles will grow to the right leaving typically
only $o(N)$ particles at each site, thus it will become harder for
particles to penetrate deep beyond the arm to completely fill a
site with $N$ particles. Thus the time required before there is
one site with $N$ particles is at least exponential in $N$. This
is the main result of this note.

\section{Formal definition of the process}

We will inductively define the following random variables.  The
variables $\{w(j,t)\}_{j,t \in \N}$ give the location of the $j$th
particle after it has take $t$ steps.  The random variables
$\{w(j,\infty)\}_{j \in \N}$ indicate where the $j$th particle
freezes. And the variables $\{f(k,i)\}_{k,i \in \N}$ indicate the
number of particles labeled  less than or equal to $i$ which have
frozen at position $k$.

We define these variables as follows.  Let $w(1,t)=0$ for all $t$.
For a fixed $j$ assume we have defined the $w(j,t)$ in such a way
that for all $j'\leq j$ we have that $\lim_{t \to \infty}w(j,t)$
exists and that we have defined $f(k,i)$ for all $i<j$. Then we
define
$$w(j,\infty)=\lim_{t \to \infty}w(j,t)$$
and
$$f(k,i) =|\{i' \leq i:\ w(i',\infty)=k\}|.$$

For any $j> 1$ let $w(j,1)=j+2$ and
$w(j,2)=j+1$.

For any $j>1$ and $t>2$ if $w(j,t)=w(j,t-1)$ then define $w(j,t+1)=w(j,t)$. If $w(j,t)\neq
w(j,t-1)$ then set
$$w(j,t+1)=\left\{%
\begin{array}{ll}
    w(j,t)+1,   & \hbox{with probability $\frac{n-f(w(j,t)-1,j-1)}{2n}$;} \\
    w(j,t)-1    & \hbox{with probability $\frac{n-f(w(j,t)-1,j-1)}{2n}$;} \\
    w(j,t),     & \hbox{with probability $\frac{f(w(j,t)-1,j-1)}{n}$;}
\\\end{array}%
\right.
$$
As simple random walk on $\Z$ is recurrent we have that $w(j,t)=w(j,t+1)$
for some $t$ almost surely.  Thus all of the random variables are well defined almost surely.

If $w(j,\infty)=k$ then we say that particle $j$ {\bf freezes} at
$k$. If for some $j$ there exists a $t$ with
$$k=w(j,t)=w(j,t+1)<\min_{t'< t}w(j,t)$$
then $w(j,\infty)=k$ and we
say that particle $j$ {\bf freezes upon arrival} at $k$. If for
some $k$ there exists $t$ with $f(k,t)=n$ then we say that there
is a {\bf blockage} at $k$.  Define the random variable
$$B=\inf \{k: \text{ there exists a blockage at $k$}\}.$$
Thus $B$ indicates the position of the leftmost blockage.

\begin{theorem} \label{main}
There exist $c>0$ such that
$$\prob(B<e^{cn})<e^{-cn}.$$
\end{theorem}

\medskip
We make the following comments about our theorem.
\begin{itemize}

\item
The proof given below clearly works for directed random walks or,
more generally, any nearest neighbor process on $\Z$. We only need
to modify the model so that the particle either freezes at some
location or disappears off to infinity.
\medskip

\item
An easy upper bound on $B$ is that for any $\epsilon>0$ there
exists $c$ such that $\prob(B>cN^N)<\epsilon$. It is of interest
to get the exact order of $B$.
\medskip

\item
The question from \cite{BY} regarding clogging of a cylinder $G
\times \Z$ is more complicated due to the geometry of the possible
cuts.  Finding the distribution on the location of the leftmost
clogging in the cylinder is an interesting open question.
\end{itemize}

\section{Proof}

For any $k$ let  $S(k)$ be the event that there exists $i$ with
$f(k,i)=n.$

\begin{lemma} \label{only}
There exists $c>0$ such that for any $k$ and $n$
 $$\prob(S(k))\leq e^{-cn}.$$
\end{lemma}
\begin{proof}
Fix $k$. For any $k',J$ and $J'$ let $C(k',J,J')$ be the number of
$j$ such that $J<j \leq J'$ and $\min_{t}w(j,t)\leq k'$.
 Let $\I$ ($I_{3n/4}, I_{7n/8}, I_{15n/16}$) be the
minimal value such that $f(k,\I)=n/2$ ($3n/4, 7n/8, 15n/16$ respectively) if such
a value exists. If $\I$ (or any of $I_{3n/4}, I_{7n/8}$, and $I_{15n/16}$)
is undefined then $S(k)$ does not occur.

Finally let
$$g(i)=\bigg|\big\{j\in (\I,i): \text{ particle $j$ freezes upon arrival at $k+1$} \big\}\bigg|.$$

Note that for all $i> \I$
$$f(k,i)) \leq n/2 + C(k,\I,i) $$
and
$$f(k+1,i)) \geq g(i)  .$$
Thus if the event $S(k)$ occurs then there must be some $i$ such
that $C(k,\I,i)\geq n/2$ and $g(i)<n$.

If $S(k)$ occurs then one of the following three things must happen:
\begin{enumerate}
  \item $g(I_{3n/4})<.22n$,
  \item $g(I_{7n/8})<.57n$, or
  \item $g(I_{15n/16})<n$.
\end{enumerate}

For every $j \in [\I,I_{3n/4}]$ with $\min_t {w(j,t)} \leq k+1$
the probability that the particle freezes upon arrival at $k+1$ is
at least one half.  Thus if the first event occurs of the first
$.47n$ particles to arrive at $k+1$ less than forty seven percent
of them freeze upon arrival at $k+1$. The probability of this is
decreasing exponentially in $n$.

For every $j \in [I_{3n/4},I_{7n/8}]$ with $\min_t {w(j,t)} \leq
k+1$ the probability that the particle freezes upon arrival at
$k+1$ is at least .75.  Thus if the first event does not occur but
the second event does, then of the first $.475n$ particles to
arrive at $k+1$ after $I_{3n/4}$ less than 74 per cent of them
freeze upon arrival at $k+1$. The conditional probability of this
given the compliment of the first event is decreasing
exponentially in $n$.

For every $j \in [I_{7n/8},I_{15n/16}]$ with $\min_t {w(j,t)} \leq
k+1$ the probability that the particle freezes upon arrival at
$k+1$ is at least .875.  Thus if the first event does not occur
but the second event does, then of the first $.4925n$ particles to
arrive at $k+1$ after $I_{3n/4}$ less than 87.4 per cent of them
freeze upon arrival at $k+1$. The conditional probability of this
given the compliment of the second event  is decreasing
exponentially in $n$.

Thus the probability that there exists $i$ such
that $C(k,\I,i)\geq n/2$ and $g(i)<n$ is exponentially small in $n$.
Thus the probability of $S(k)$ is as well.
\end{proof}

\begin{pfofthm}{\ref{main}}
This follows  by replacing $c$ with $c/2$  in Lemma \ref{only}.
\end{pfofthm}

\end{document}